\documentclass{amsart}
\usepackage{amsfonts}
\usepackage{amssymb}
\usepackage{times}
\usepackage{xcolor}
\usepackage{amsmath}
\usepackage{amsthm}
\usepackage{comment}
\usepackage{xcolor}
\usepackage{graphicx}

\makeatletter
\def\Ddots{\mathinner{\mkern1mu\raise\p@
\vbox{\kern7\p@\hbox{.}}\mkern2mu
\raise4\p@\hbox{.}\mkern2mu\raise7\p@\hbox{.}\mkern1mu}}
\makeatother

\newtheorem{theorem}{Theorem}[section]

\newtheorem{proposition}{Proposition}[section]
\newtheorem{definition}{Definition}[section]

\newtheorem{conjecture}{Conjecture}[section]

\newenvironment{dem1,1}[1][Proof of Theorem \ref{mainthm1,1}]{\noindent\textit{#1.} }{\hfill $\square$}
\newenvironment{dem1,2}[1][Proof of Theorem \ref{mainthm1,2}]{\noindent\textit{#1.} }{\hfill $\square$}
\newenvironment{dem2,1}[1][Proof of Theorem \ref{mainthm2,1}]{\noindent\textit{#1.} }{\hfill $\square$}
\newenvironment{dem2,2}[1][Proof of Theorem \ref{mainthm2,2}]{\noindent\textit{#1.} }{\hfill $\square$}

\begin{document}

\title{Unraveling Functional Equations in Composition Algebra: Resolving Conjectures and Examining Implications}

\thanks{The first author was supported by FAPESP number 2022/14579-0.}

\author[Daniel Kawai]{Daniel Kawai}
\address{Daniel Kawai, S\~{a}o Paulo University, 
Rua Mat\~{a}o, 1010, 
CEP 05508-090, S\~{a}o Paulo, Brazil\\
{\em E-mail}: {\tt daniel.kawai@usp.br}}{}

\author[Bruno Leonardo Macedo Ferreira]{Bruno Leonardo Macedo Ferreira}
\address{Bruno Leonardo Macedo Ferreira, Federal University of Technology-Paran\'{a}-Brazil, Federal University of ABC and S\~{a}o Paulo University, 
Avenida Professora Laura Pacheco Bastos, 800, 
85053-510, Guarapuava, Brazil, Av. dos Estados, 5001 - Bangú, Santo André - SP, 09280-560 and Rua Mat\~{a}o, 1010, 
CEP 05508-090, S\~{a}o Paulo, Brazil \\\linebreak
{\em E-mail}: {\tt 
 brunolmfalg@gmail.com}}{}

\begin{abstract}
We address the conjectures left by the recent article by Ferreira et al. titled ``Commuting maps and identities with inverses on alternative division rings.'' We also present an example showing the necessity of the conditions of the results that answer the conjectures.
\vspace*{.1cm}

\vspace*{.1cm}

\end{abstract}

\maketitle

\section{Introduction}
In 1987, Vukman \cite{Vuk} studied a result about Cauchy functional equations in associative division rings, and he showed that if $D$ is a noncommutative associative division ring with characteristics different from $2$ and $f : D \rightarrow D$ is an additive function satisfying:
\begin{equation}
\label{vukman}
f(x) = -x^2f(x^{-1})
\end{equation}
for every $x \in D^\times$, where $D^\times$ is the set of invertible elements of $D$, then $f(x) = 0$ for all $x \in D$.

In 2024, Ferreira et al. \cite{smigly} studied the functional equation \eqref{vukman} in alternative division rings and proved the following.

\begin{theorem}\label{th5.1}
Let $D$ be a noncommutative alternative division ring with characteristic different from $2$, and let $f : D \rightarrow D$ be an additive map satisfying the identity \eqref{vukman}:
\begin{equation*}
f(x) = -x^2f(x^{-1})
\end{equation*}
for for every $x \in D^\times$. Then $f(x) = 0$ for all $x \in D$.
\end{theorem}

They also posed the following conjectures:

\begin{conjecture}\label{c1}
Is Theorem \ref{th5.1} true if $D = O$ is a split-octonion algebra?
\end{conjecture}

\begin{conjecture}\label{c2}
Let $D$ be an alternative division ring and suppose $f, g : D \rightarrow D$ are additive maps satisfying:
\begin{equation}
\label{vukmana}
f(x) + x^2g(x^{-1}) = 0
\end{equation}
for all $x \in D^\times$. Then there exists a fixed element $q \in D$ such that $f(x) = xq$
and $g(x) = -xq$ for all $x \in D$.
\end{conjecture}

Note that the articles \cite{Vuk} and \cite{smigly} deal with the particular case of equation \eqref{vukmana} where $f=g$. The central focus of this work is to provide answers to these conjectures. Our approach aims to characterize the additive mappings that preserve that identity, exploring both theoretical aspects and practical implications.

We hope that this study will contribute not only to the understanding of functional equation theory but also to the advancement of knowledge in various areas of mathematics.

\section{Preliminaries}

Let $R$ be a ring (not necessarily associative or commutative) and consider the following notational convention for its multiplication operation: that $xy \cdot z = (xy)z$ and $x\cdot yz = x(yz)$ for $x, y, z \in R$, simply to reduce the number of necessary parentheses. We denote the associator of $R$ by $(x, y, z) = xy \cdot z - x \cdot yz$ for $x, y, z \in R$.

A ring $R$ is said to be \emph{alternative} if $(x, x, y) = 0 = (y, x, x)$ for all $x, y \in R$. One easily sees that any associative ring is an alternative ring. It is also well-known that every alternative ring satisfies the flexible identity $(x, y, x) = 0$, for all $x, y \in R$. Alternative rings also satisfy the Moufang identities
\begin{equation*}
(xax)y = x[a(xy)],
\end{equation*}
\begin{equation*}
y(xax) = [(yx)a]x,
\end{equation*}
\begin{equation*}
(xy)(ax) = x(ya)x
\end{equation*}
for all $x, y, a \in R$. If parentheses or the symbol $\cdot$ are not needed to write an expression, we omit them.

In any ring $R$ with $1$, an element $x$ is said to have an inverse $x^{-1}$ if there is an element $x^{-1}$ in $R$ satisfying $xx^{-1} = x^{-1}x = 1$. The set of all invertible elements of $R$ is denoted $R^\times$. Just as in the associative case, inverses are unique in alternative rings. If every nonzero element in an alternative ring $R$ with $1$ has an inverse, then $R$ is called an alternative division ring. If $R$ is an alternative ring with $1$, then a subset $S\subseteq R$ is said to be a \emph{division subring} if $S$ is a subring of $R$, $1\in S$ and for every $x\in S\cap R^\times$ we have $x^{-1}\in S$.

It is well known that, if $R$ is an alternative division ring, then for every $x,y\in R$ the division subring generated by $\{x,y\}$ is associative, as the following proposition.

\begin{proposition}
In an alternative ring with $1$, the division subring generated by any two elements is associative.
\end{proposition}
\begin{proof}
Let $A$ be an alternative ring with $1$. We will show that if $(a,b,c)=0$ and $a\in R^\times$, then $(a^{-1},b,c)=0$. Indeed, we have:
\[
(a^{-1},b,c)=(a^{-1},a^{-1}(ab),c)=(a^{-1},ab,c)a^{-1},
\]
but we have:
\[
(a^{-1},ab,c)=-(b,aa^{-1},c)+b(a^{-1},a,c)+a^{-1}(b,a,c)=0,
\]
so $(a^{-1},b,c)=0$.

Let $u,v\in A$ be two elements. Let $R_k$ be the submodule generated by non-associative words of the form $w(a_1,\dots,a_r)$ where there exist $k_1,\dots,k_r\geq 1$ such that $k_1+\cdots+k_r\leq k$ and $a_i\in\{u,v\}$ if $k_i=1$ and $a_i\in (R_{k_i-1}\cap R^\times)^{-1}$ if $k_i\geq 2$. Let $S=\bigcup_{k=1}^\infty R_k$. Then $R_1\subseteq R_2\subseteq\cdots$ and also $S$ is a submodule. For $x\in R_k$ and $y\in R_l$, then $x$ is a linear combination of $w_x(a_1,\dots,a_r)$ where there exist $k_1,\dots,k_s,l_1,\dots,l_s\geq1$ such that $k_1+\cdots+k_r+l_1+\cdots+l_s\leq k+l$ and $a_i\in\{u,v\}$ if $k_i=1$ and $a_i\in (R_{k_i-1}\cap R^\times)^{-1}$ if $k_i\geq 2$, and similarly for $y_j$. Therefore $xy$ is a linear combination of words $w_x(a_1,\dots,a_r)w_y(b_1,\dots,b_s)$ where there exist $k_1,\dots,k_s,l_1,\dots,l_s\geq1$ such that $k_1+\cdots+k_r+l_1+\cdots+l_s\leq k+l$ and $a_i\in\{u,v\}$ if $k_i=1$ and $a_i\in (R_{k_i-1}\cap R^\times)^{-1}$ if $k_i\geq 2$, and similarly for $b_j$. Hence $xy\in R_{k+l}$. Finally, if $x\in R_k\cap R^\times$, then $x^{-1}\in (R_k\cap R^\times)^{-1}$, so $x^{-1}\in R_{k+1}$. Also, it is easy to see that $u,v\in R_1$. Therefore $S$ is the division subring generated by $u$ and $v$.

We will show by induction on $d(x)+d(y)+d(z)$ that $(x,y,z)=0$ for any $x,y,z\in S$, where $d(x)$ is the smallest $k$ such that $x\in R_k$. If $d(x)+d(y)+d(z)=3$, then $d(x)=d(y)=d(z)=1$, so $x$, $y$, and $z$ are linear combinations of $u$ and $v$, so $(x,y,z)=0$. Now if $d(x)+d(y)+d(z)>3$, we can assume that $x$, $y$, and $z$ are words $w_x(a_1,\dots,a_r)$, $w_y(b_1,\dots,b_s)$, and $w_z(c_1,\dots,c_t)$ and that $x_1,\dots,x_r,y_1,\dots,y_s,z_1,\dots,z_t\geq1$ satisfy $x_1+\dots+x_r\leq d(x)$ and $a_i\in\{u,v\}$ if $x_i=1$ and $a_i\in (R_{x_i-1}\cap R^\times)^{-1}$ if $x_i\geq 2$, and similarly for $y_j$ and $z_k$. Since $n=d(x)+d(y)+d(z)$, then $d(x),d(y),d(z)<n$, so by the induction hypothesis, we can assume that $w_x(a_1,\dots,a_r)=(\cdots(a_1a_2)a_3\cdots)a_r$ and so on. We have two cases:

\textbf{Case 1:} If $x_r\geq 2$, then $a_r=a^{-1}$ for some $a\in R_{x_r-1}\cap R^\times$. 
\begin{enumerate}
    \item If $r=1$, then $(x,y,z)=(a^{-1},y,z)=0$, because $d(a)\leq x_r-1\leq d(x)-1$ and so $(a,y,z)=0$.
    \item If $r\geq 2$, let $\overline{a}=(\cdots(a_1a_2)a_3\cdots)a_{r-1}$, then:
    \[
    \begin{aligned}
    (x,y,z)&=(\overline{a}a,y,z)\\
    &=-(\overline{a}y,a_r,z)+a_r(\overline{a},y,z)+y(a,\overline{a},z),
    \end{aligned}
    \]
    but $d(\overline{a})\leq d(a_1)+\cdots+d(a_{r-1})\leq x_1+\cdots+x_{r-1}$ and $d(a_r)\leq x_r$, so $d(\overline{a})<d(x)$, hence:
    \[
    \begin{aligned}
    (x,y,z)&=-(\overline{a}y,a_r,z)\\
    &=-u(\overline{a},\overline{b}z,u)=0,
    \end{aligned}
    \]
    since $(\overline{a}y,a,z)=0$ as $d(a)\leq x_r-1<x_r$ so $d(\overline{a}y)+d(a)+d(z)<n$.
\end{enumerate}

\textbf{Case 2:} If $x_r=y_s=z_t=1$, then $a_r,b_s,c_t\in\{u,v\}$, so we can assume without loss of generality that $a_r=b_s=u$, then:
\[
\begin{aligned}
(x,y,z)&=(\overline{a}u,\overline{b}u,z)\\
&=-(\overline{a}u,\overline{b}z,u)+u(\overline{a}u,\overline{b},z)+z(\overline{a}u,\overline{b},u)\\
&=-(\overline{a}u,\overline{b}z,u)\\
&=-u(\overline{a},\overline{b}z,u)=0.
\end{aligned}
\]
\end{proof}

It is also pertinent to revisit the well-known identity attributed to Hua, an important element for the development of this paper. Let $R$ be an unital alternative ring and let $a$, $b$ be any two elements of $R$ with $a,b,ab-1\in R^\times$. Then Hua’s identity holds:
\begin{equation*}
a - (a^{-1} + (b^{-1} - a)^{-1})^{-1} = aba.
\end{equation*}

\section{Alterntative Division Rings}
In this section, we present theorems that address the Conjecture \ref{c2}.

\begin{theorem}\label{t1}
Let $D$ satisfy one of the following properties:
\begin{itemize}
\item[a)] $D$ is a noncommutative alternative division ring with $\mathrm{Char}(D)\neq2$.
\item[b)] $D$ is a perfect field with $\mathrm{Char}(D)=2$.
\end{itemize}
Suppose $f, g : D \rightarrow D$ are additive maps satisfying \eqref{vukmana}:
\begin{equation*}
f(x) + x^2g(x^{-1}) = 0
\end{equation*}
for all $x \in D^\times$. Then there exists a fixed element $q \in D$ such that $f(x) = xq$ and $g(x) = -xq$ for all $x \in D$.
\end{theorem}

\begin{proof}
We proceed with the proof of Theorem \ref{t1} by delineating several steps.

\medskip
\textbf{Step 1:}
Consider an unital alternative ring $D$, and let $f, g : D \rightarrow D$ be additive maps satisfying \eqref{vukmana} for all $x \in D^\times$. Choose $a\in D$ such that $a,a-1\in D^\times$. By invoking Hua's identity and exploiting the associativity of the division subring generated by $\{a,x\}$ for every $x\in D$, we derive:
\begin{equation}
a^2=a-(a^{-1}+(1-a)^{-1})^{-1}
\end{equation}
and consequently:
\begin{equation*}
\begin{aligned}
f(a^2)&=f(a-(a^{-1}+(1-a)^{-1})^{-1})\\
&=f(a)-f((a^{-1}+(1-a)^{-1})^{-1})\\
&\overset{\eqref{vukmana}}{=}f(a)+(a^{-1}+(1-a)^{-1})^{-2}g(a^{-1}+(1-a)^{-1})\\
&=f(a)+a^2(1-a)^2g(a^{-1}+(1-a)^{-1})\\
&=f(a)+a^2(1-a)^2(g(a^{-1})+g((1-a)^{-1}))\\
&=f(a)+a^2(1-a)^2g(a^{-1})+a^2(1-a)^2g((1-a)^{-1})\\
&\overset{\eqref{vukmana}}{=}f(a)-a^2(1-a)^2a^{-2}f(a)-a^2(1-a)^2(1-a)^{-2}f(1-a)\\
&=f(a)-(1-a)^2f(a)-a^2f(1-a)\\
&=f(a)-(1-a)^2f(a)-a^2(f(1)-f(a))\\
&=f(a)-(1-a)^2f(a)-a^2f(1)+a^2f(a)\\
&=(1+a^2-(1-a)^2)f(a)-a^2f(1)\\
&=2af(a)-a^2f(1),
\end{aligned}
\end{equation*}
which can be succinctly expressed as:
\begin{equation}
\label{4}
f(a^2)=2af(a)-a^2f(1).
\end{equation}

\textbf{Step 2:} Assume $D$ satisfies either (a) or (b), and let $f, g : D \rightarrow D$ be additive maps satisfying \eqref{vukmana} for all $x \in D^\times$. In this scenario, $D$ qualifies as an alternative division ring. Consequently, for any $a\in D$, if $a\neq0$ and $a\neq 1$, then $a\in D^\times$ and $a-1\in D^\times$, allowing us to apply \eqref{4}. Moreover, \eqref{4} remains valid for $a=0$ or $a=1$. Define $h(x)=f(x)-xf(1)$ for all $x\in D$. Then $h:D\rightarrow D$ is additive and satisfies:
\begin{equation}
\label{5a}
h(x^2)=2xh(x)
\end{equation}
for every $x\in D$.
\begin{itemize}
\item[a)] In case (a), leveraging the Bruck and Kleinfeld Theorem, $D$ manifests as either an associative division algebra or an octonion division algebra over its center. We can adopt the proof strategy from the Main Theorem of \cite{Vuk} for the associative case and Theorem 5.1 of \cite{smigly} for the octonion case to establish $h(x)=0$ for all $x\in D$.
\item[b)] In case (b), for every $x\in D$, there exists $r\in D$ such that $r^2=x$, implying $h(x)=h(r^2)=2rh(r)=0$.
\end{itemize}
Consequently, in both scenarios, we infer:
\begin{equation}
\label{6a}
f(x)=xf(1)
\end{equation}
for all $x\in D$. Thus, for every $x\neq0$, by invoking \eqref{vukmana} with $x^{-1}$ instead of $x$, we deduce:
\begin{equation*}
f(x^{-1})+x^{-2}g(x)=0,
\end{equation*}
which, in turn, yields:
\begin{equation}
\label{7a}
g(x)=-xf(1).
\end{equation}
Equation \eqref{7a} holds for $x=0$ as well.
\end{proof}

Now, we will prove the necessity of the condition of perfectness in Theorem \ref{t1}, case (b), by illustrating that in the non-perfect field $D = \mathbb{Z}_2(t)$ of rational functions with coefficients in $\mathbb{Z}_2$, there exist numerous solutions $(f, g)$ to the functional equation \eqref{vukmana}, where $f$ and $g$ are additive, aside from the standard one.

\begin{theorem}
Let $D=\mathbb{Z}_2(t)$ be the field of rational functions with coefficients in $\mathbb{Z}_2$. For any $A,B\in D$, there is exactly one pair $(f,g)$ consisting of additive functions $f,g:D\rightarrow D$ that satisfy $f(1)=A$, $f(t)=B$ and the identity \eqref{vukmana}:
\begin{equation*}
f(x) + x^2g(x^{-1}) = 0
\end{equation*}
for every $x\in D$ such that $x\neq0$.
\end{theorem}
\begin{proof}
\textbf{1) Uniqueness: }
Let $f,g:D\rightarrow D$ be additive functions that satisfy $f(1)=A$, $f(t)=B$ and \eqref{vukmana} for all $x\neq0$. Then, for every $a,b\in D$ such that $ab\neq0,1$, by Hua identity and the fact that $D$ is a field with $\mathrm{Char}(D)=2$:
\begin{equation*}
a^2b=a+(a^{-1}+(b^{-1}+a)^{-1})^{-1}
\end{equation*}
we have:
\begin{equation*}
\begin{array}{rcl}
f(a^2b)
&=&f(a+(a^{-1}+(b^{-1}+a)^{-1})^{-1})\\
&=&f(a)+f((a^{-1}+(b^{-1}+a)^{-1})^{-1})\\
&\overset{\eqref{vukmana}}{=}&f(a)+(a^{-1}+(b^{-1}+a)^{-1})^{-2}g(a^{-1}+(b^{-1}+a)^{-1})\\
&=&f(a)+(a^2b+a)^2g(a^{-1}+(b^{-1}+a)^{-1})\\
&=&f(a)+(a^2b+a)^2(g(a^{-1})+g((b^{-1}+a)^{-1}))\\
&=&f(a)+(a^2b+a)^2g(a^{-1})+(a^2b+a)^2g((b^{-1}+a)^{-1})\\
&\overset{\eqref{vukmana}}{=}&f(a)+(a^2b+a)^2a^{-2}f(a)+(a^2b+a)^2(b^{-1}+a)^{-2}f(b^{-1}+a)\\
&=&f(a)+(ab+1)^2f(a)+a^2b^2(a+b^{-1})^2(b^{-1}+a)^{-2}f(b^{-1}+a)\\
&=&f(a)+(a^2b^2+1)f(a)+a^2b^2f(b^{-1}+a)\\
&=&a^2b^2f(a)+a^2b^2(f(b^{-1})+f(a))\\
&=&a^2b^2f(a)+a^2b^2f(b^{-1})+a^2b^2f(a)\\
&\overset{\eqref{vukmana}}{=}&a^2b^2f(a)+a^2b^2b^{-2}g(b)+a^2b^2f(a)\\
&=&a^2b^2f(a)+a^2g(b)+a^2b^2f(a)
\\
&=&a^2g(b),
\end{array}
\end{equation*}
that is:
\begin{equation}
\label{14}
f(a^2b)=a^2g(b).
\end{equation}
By $f(0)=g(0)=0$ and \eqref{vukmana}, the above equation \eqref{14} is also true for $a,b\in D$ such that $ab=0$ or $ab=1$. So \eqref{14} holds for all $a,b\in D$. By \eqref{14} with $a=1$, for all $b\in D$ we have $f(b)=g(b)$. So $f=g$ and for all $a\in D$ we have:
\begin{equation}
\label{15}
f(a^2b)=a^2f(b).
\end{equation}
For every $n\in\mathbb{Z}$, by \eqref{15} with $a=t$ and $b=t^n$, we have:
\begin{equation}
\label{16}
f(t^{n+2})=t^2f(t^n).
\end{equation}
Therefore, using induction in $n$ and \eqref{16}, for all $n\in\mathbb{Z}$ we have:
\begin{equation*}
f(t^{2n})=t^{2n}A,\quad\quad f(t^{2n+1})=t^{2n}B.
\end{equation*}
Also, if $P(t)$ is a polynomial with coefficients in $\mathbb{Z}_2$ and $p\in\mathbb{Z}$, then it is easy to see that:
\begin{equation}
\label{18}
f(P(t^2)t^{2p})=P(t^2)t^{2p}A,\quad\quad f(P(t^2)t^{2p+1})=P(t^2)t^{2p}B
\end{equation}
and also:
\begin{equation}
\label{d3}
P(t^2)=P(t)^2.
\end{equation}
Every element of $x\in D$ can be represented by the form:
\begin{equation*}
x=\frac{P(t^2)+Q(t^2)t}{R(t^2)+S(t^2)t},
\end{equation*}
where $P(t),Q(t),R(t),S(t)$ are polynomials with coefficients in $\mathbb{Z}_2$, so, by representing:
\begin{equation*}
P(t)=\sum_{p=0}^{\hat{p}}P_pt^p,\quad\quad Q(t)=\sum_{q=0}^{\hat{q}}Q_qt^q,
\end{equation*}
where $P_p,Q_q\in\mathbb{Z}_2$ we have:
\begin{equation*}
\begin{aligned}
f(x)&=f\left(\frac{P(t^2)+Q(t^2)t}{R(t^2)+S(t^2)t}\right)\\
&=\sum\limits_{p=0}^{\hat{p}}P_pf\left(\frac{t^{2p}}{R(t^2)+S(t^2)t}\right)+\sum\limits_{q=0}^{\hat{q}}Q_qf\left(\frac{t^{2q+1}}{R(t^2)+S(t^2)t}\right)\\
&\overset{\text{\eqref{vukmana}}}{=}\sum\limits_{p=0}^{\hat{p}}\frac{P_pt^{4p}}{(R(t^2)+S(t^2)t)^2}f\left(\frac{R(t^2)+S(t^2)t}{t^{2p}}\right)+\sum\limits_{q=0}^{\hat{q}}\frac{Q_qt^{4q+2}}{(R(t^2)+S(t^2)t)^2}f\left(\frac{R(t^2)+S(t^2)t}{t^{2q+1}}\right)\\
&=\sum\limits_{p=0}^{\hat{p}}\frac{P_pt^{4p}f(R(t^2)t^{-2p}+S(t^2)t^{-2p+1})}{(R(t^2)+S(t^2)t)^2}+\sum\limits_{q=0}^{\hat{q}}\frac{Q_qt^{4q+2}f(R(t^2)t^{-2q-1}+S(t^2)t^{-2q})}{(R(t^2)+S(t^2)t)^2}\\
&\overset{\text{\eqref{18}}}{=}\sum\limits_{p=0}^{\hat{p}}\frac{P_pt^{4p}(R(t^2)t^{-2p}A+S(t^2)t^{-2p}B)}{(R(t^2)+S(t^2)t)^2}+\sum\limits_{q=0}^{\hat{q}}\frac{Q_qt^{4q+2}(R(t^2)t^{-2q-2}B+S(t^2)t^{-2q}A)}{(R(t^2)+S(t^2)t)^2}\\
&=\sum\limits_{p=0}^{\hat{p}}\frac{P_pt^{2p}(R(t^2)A+S(t^2)B)}{(R(t^2)+S(t^2)t)^2}+\sum\limits_{q=0}^{\hat{q}}\frac{Q_qt^{2q}(R(t^2)B+S(t^2)t^2A)}{R(t^2)+S(t^2)t}\\
&=\frac{P(t^2)(R(t^2)A+S(t^2)B)}{(R(t^2)+S(t^2)t)^2}+\frac{Q(t^2)(R(t^2)B+S(t^2)t^2A)}{(R(t^2)+S(t^2)t)^2}\\
&=\frac{P(t^2)R(t^2)+Q(t^2)S(t^2)t^2}{(R(t^2)+S(t^2)t)^2}A+\frac{P(t^2)S(t^2)+Q(t^2)R(t^2)}{(R(t^2)+S(t^2)t)^2}B\\
&=\left(\frac{P(t)R(t)+Q(t)S(t)t}{R(t^2)+S(t^2)t}\right)^2A+\left(\frac{P(t)S(t)+Q(t)R(t)}{R(t^2)+S(t^2)t}\right)^2B,
\end{aligned}
\end{equation*}
therefore:
\begin{equation}
\label{17}
\begin{array}{rcl}
f\left(\frac{P(t^2)+Q(t^2)t}{R(t^2)+S(t^2)t}\right)&=&\left(\frac{P(t)R(t)+Q(t)S(t)t}{R(t^2)+S(t^2)t}\right)^2A+\left(\frac{P(t)S(t)+Q(t)R(t)}{R(t^2)+S(t^2)t}\right)^2B
\end{array}
\end{equation}
\textbf{2) Existence: }
It remains now to show that for any $A,B\in D$ there exists a function $f:D\rightarrow D$ such that \eqref{17} holds for any polynomials $P(t),Q(t),R(t),S(t)$ with coefficients in $\mathbb{Z}_2$, and to show that in fact $f$ is additive, $f(1)=A$, $f(t)=B$ and, if $g=f$, then $f$ and $g$ satisfy \eqref{vukmana} for all $x\neq0$. We will divide this task into steps.

\medskip
\textbf{2,1)} We will show that there exists a function $f:D\rightarrow D$ satisfying \eqref{17} for  any polynomials $P(t),Q(t),R(t),S(t)$ with coefficients in $\mathbb{Z}_2$. In other words, we will show that the function $f$ as described in \eqref{17} is well defined. 

\medskip
Let $P(t),\hat{P}(t),Q(t),\hat{Q}(r),R(t),\hat{R}(t),S(t),\hat{S}(t)$ be polynomials such that:
\begin{equation*}
\frac{P(t^2)+Q(t^2)t}{R(t^2)+S(t^2)t}=\frac{\hat{P}(t^2)+\hat{Q}(t^2)t}{\hat{R}(t^2)+\hat{S}(t^2)t}
\end{equation*}
Then:
\begin{equation*}
(P(t^2)+Q(t^2)t)(\hat{R}(t^2)+\hat{S}(t^2)t)=(\hat{P}(t^2)+\hat{Q}(t^2)t)(R(t^2)+S(t^2)t)
\end{equation*}
Hence:
\begin{equation*}
(P(t^2)\hat{R}(t^2)+Q(t^2)\hat{S}(t^2)t^2)+(P(t^2)\hat{S}(t^2)+Q(t^2)\hat{R}(t^2))t
\end{equation*}
\begin{equation*}
=(\hat{P}(t^2)R(t^2)+\hat{Q}(t^2)S(t^2)t^2)+(\hat{P}(t^2)S(t^2)+\hat{Q}(t^2)R(t^2))t
\end{equation*}
Both sides of the above equality are polynomials over $t$, so we can compare their coefficients and obtain the following equalities:
\begin{equation}
\label{d1}
P(t)\hat{R}(t)+Q(t)\hat{S}(t)t=\hat{P}(t)R(t)+\hat{Q}(t)S(t)t
\end{equation}
\begin{equation}
\label{d2}
P(t)\hat{S}(t)+Q(t)\hat{R}(t)=\hat{P}(t)S(t)+\hat{Q}(t)R(t)
\end{equation}
Multiplying \eqref{d1} by $R(t)\hat{R}(t)$ and by $S(t)\hat{S}(t)t$, and multiplying \eqref{d2} by $R(t)\hat{S}(t)t$ and by $S(t)\hat{R}(t)t$, we obtain respectively:
\begin{equation*}
P(t)R(t)\hat{R}(t)^2+Q(t)R(t)\hat{R}(t)\hat{S}(t)t=\hat{P}(t)\hat{R}(t)R(t)^2+\hat{Q}(t)\hat{R}(t)R(t)S(t)t
\end{equation*}
\begin{equation*}
P(t)S(t)\hat{R}(t)\hat{S}(t)t+Q(t)S(t)\hat{S}(t)^2t^2=\hat{P}(t)\hat{S}(t)R(t)S(t)t+\hat{Q}(t)\hat{S}(t)S(t)^2t^2
\end{equation*}
\begin{equation*}
P(t)R(t)\hat{S}(t)^2t+Q(t)R(t)\hat{R}(t)\hat{S}(t)t=\hat{P}(t)\hat{S}(t)R(t)S(t)t+\hat{Q}(t)\hat{S}(t)R(t)^2t
\end{equation*}
\begin{equation*}
P(t)S(t)\hat{R}(t)\hat{S}(t)t+Q(t)S(t)\hat{R}(t)^2t=\hat{P}(t)\hat{R}(t)S(t)^2t+\hat{Q}(t)\hat{R}(t)R(t)S(t)t
\end{equation*}
By summing up the above four equations, and using $\mathrm{Char}(D)=2$ and \eqref{d3}, we obtain:
\begin{equation*}
(P(t)R(t)+Q(t)S(t)t)(\hat{R}(t^2)+\hat{S}(t^2)t)=(\hat{P}(t)\hat{R}(t)+\hat{Q}(t)\hat{S}(t)t)(R(t^2)+S(t^2)t)
\end{equation*}
so that:
\begin{equation}
\label{ea1}
\frac{P(t)R(t)+Q(t)S(t)t}{R(t^2)+S(t^2)t}=\frac{\hat{P}(t)\hat{R}(t)+\hat{Q}(t)\hat{S}(t)t}{\hat{R}(t^2)+\hat{S}(t^2)t}
\end{equation}
Now, multiplying \eqref{d1} by $S(t)\hat{R}(t)$ and by $R(t)\hat{S}(t)$, and multiplying \eqref{d2} by $S(t)\hat{S}(t)t$ and by $R(t)\hat{R}(t)$, we obtain respectively:
\begin{equation*}
P(t)S(t)\hat{R}(t)^2+Q(t)S(t)\hat{R}(t)\hat{S}(t)t=\hat{P}(t)\hat{R}(t)R(t)S(t)+\hat{Q}(t)\hat{R}(t)S(t)^2t
\end{equation*}
\begin{equation*}
P(t)R(t)\hat{R}(t)\hat{S}(t)+Q(t)R(t)\hat{S}(t)^2t=\hat{P}(t)\hat{S}(t)R(t)^2+\hat{Q}(t)\hat{S}(t)R(t)S(t)t
\end{equation*}
\begin{equation*}
P(t)S(t)\hat{S}(t)^2t+Q(t)S(t)\hat{R}(t)\hat{S}(t)t=\hat{P}(t)\hat{S}(t)S(t)^2t+\hat{Q}(t)\hat{S}(t)R(t)S(t)t
\end{equation*}
\begin{equation*}
P(t)R(t)\hat{R}(t)\hat{S}(t)+Q(t)R(t)\hat{R}(t)^2=\hat{P}(t)\hat{R}(t)R(t)S(t)+\hat{Q}(t)\hat{R}(t)R(t)^2
\end{equation*}
By adding the above four equations, and using $\mathrm{Char}(D)=2$ and \eqref{d3}, we have:
\begin{equation*}
(P(t)S(t)+Q(t)R(t))(\hat{R}(t^2)+\hat{S}(t^2)t)=(\hat{P}(t)\hat{S}(t)+\hat{Q}(t)\hat{R}(t))(R(t^2)+S(t^2)t)
\end{equation*}
so that:
\begin{equation}
\label{ea2}
\frac{P(t)S(t)+Q(t)R(t)}{R(t^2)+S(t^2)t}=\frac{\hat{P}(t)\hat{S}(t)+\hat{Q}(t)\hat{R}(t)}{\hat{R}(t^2)+\hat{S}(t^2)t}
\end{equation}
We now square both sides of \eqref{ea1} and multiply by $A$, square both sides of \eqref{ea2} and multiply by $B$, and add the equations thus obtained, to obtain:
\begin{equation*}
\left(\frac{P(t)R(t)+Q(t)S(t)t}{R(t^2)+S(t^2)t}\right)^2A+\left(\frac{P(t)S(t)+Q(t)R(t)}{R(t^2)+S(t^2)t}\right)^2B
\end{equation*}
\begin{equation*}
=\left(\frac{\hat{P}(t)\hat{R}(t)+\hat{Q}(t)\hat{S}(t)t}{\hat{R}(t^2)+\hat{S}(t^2)t}\right)^2A+\left(\frac{\hat{P}(t)\hat{S}(t)+\hat{Q}(t)\hat{R}(t)}{\hat{R}(t^2)+\hat{S}(t^2)t}\right)^2B
\end{equation*}
\textbf{2,2)} Let us show that $f$ is additive. Let $x,y\in D$, then, by bringing the fractions to a common denominator, we can have polynomials $P(t),\hat{P}(t), Q(t),\hat{Q}(t), R(t), S(t)$ such that:
\begin{equation*}
x=\frac{P(t^2)+Q(t^2)t}{R(t^2)+S(t^2)t},\quad\quad y=\frac{\hat{P}(t^2)+\hat{Q}(t^2)t}{R(t^2)+S(t^2)t},
\end{equation*}
thus:
\begin{equation*}
x+y=\frac{(P(t^2)+\hat{P}(t^2))+(Q(t^2)+\hat{Q}(t^2))t}{R(t^2)+S(t^2)t}
\end{equation*}
therefore, using the fact that $x\mapsto x^2$ is additive in fields of characteristic $2$:
\begin{equation*}
\begin{aligned}
f(x+y)&\overset{\eqref{17}}{=}\left(\frac{(P(t)+\hat{P}(t))R(t)+(Q(t)+\hat{Q}(t))S(t)t}{R(t^2)+S(t^2)t}\right)^2A\\
&+\left(\frac{(P(t)+\hat{P}(t))S(t)+(Q(t)+\hat{Q}(t))R(t)}{R(t^2)+S(t^2)t}\right)^2B\\
&=\left(\frac{P(t)R(t)+Q(t)S(t)t}{R(t^2)+S(t^2)t}\right)^2A+\left(\frac{P(t)S(t)+Q(t)R(t)}{R(t^2)+S(t^2)t}\right)^2B\\
&+\left(\frac{\hat{P}(t)R(t)+\hat{Q}(t)S(t)t}{\hat{R}(t^2)+\hat{S}(t^2)t}\right)^2A+\left(\frac{\hat{P}(t)S(t)+\hat{Q}(t)R(t)}{R(t^2)+S(t^2)t}\right)^2B\\
&\overset{\eqref{17}}{=}f(x)+f(y).
\end{aligned}
\end{equation*}
\textbf{2,3)} Let us show that $f(1)=A$ and $f(t)=B$. Indeed, we have:
\begin{equation*}
1=\frac{1+0t}{1+0t}
\end{equation*}
so that:
\begin{equation*}
f(1)=\left(\frac{1\cdot1+0\cdot0t}{1+0t}\right)^2A+\left(\frac{1\cdot 0+0\cdot 1}{1+0t}\right)^2B=A
\end{equation*}
Also, we have:
\begin{equation*}
t=\frac{0+1t}{1+0t}
\end{equation*}
thereby:
\begin{equation*}
f(t)=\left(\frac{0\cdot1+1\cdot0t}{1+0t}\right)^2A+\left(\frac{0\cdot 0+1\cdot 1}{1+0t}\right)^2B=B
\end{equation*}
\textbf{2,4)} Let us show that, if $g=f$, then $(f,g)$ satisfy \eqref{vukmana} for all $x\neq0$. Let $x\neq0$, then we have polynomials $P(t),Q(t),R(t),S(t)$ such that:
\begin{equation*}
x=\frac{P(t^2)+Q(t^2)t}{R(t^2)+S(t^2)t},
\end{equation*}
hence:
\begin{equation*}
x^{-1}=\frac{R(t^2)+S(t^2)t}{P(t^2)+Q(t^2)t},
\end{equation*}
therefore:
\begin{equation*}
\begin{aligned}
x^2f(x^{-1})&\overset{\eqref{17}}{=}\left(\frac{P(t^2)+Q(t^2)t}{R(t^2)+S(t^2)t}\right)^2\left(\left(\frac{R(t)P(t)+S(t)Q(t)t}{P(t^2)+Q(t^2)t}\right)^2A+\left(\frac{R(t)Q(t)+S(t)P(t)}{P(t^2)+R(t^2)t}\right)^2B\right)\\
&=\left(\frac{P(t)R(t)+Q(t)S(t)t}{R(t^2)+S(t^2)t}\right)^2A+\left(\frac{P(t)S(t)+Q(t)R(t)}{R(t^2)+S(t^2)t}\right)^2B\\
&\overset{\eqref{17}}{=}f(x),
\end{aligned}
\end{equation*}
so that, by $\mathrm{Char}(D)=2$, we obtain:
\begin{equation*}
f(x)+x^2f(x^{-1})=0.
\end{equation*}
Therefore, if $g=f$, then $(f,g)$ satisties \eqref{vukmana} for all $x\neq0$.

\medskip
\textbf{Conclusion}:
Thus, given any $A, B\in D$, we have proven both the uniqueness and existence of the additive functions $f,g: D\rightarrow D$ satisfying \eqref{vukmana} for every $x\neq0$. Therefore, the theorem is established.
\end{proof}

\section{Split Octonion Algebras}

In this section, we delve into investigating Conjecture \ref{c1}, which questions whether the outcome of Theorem 3.1 extends to split octonion algebras. However, before delving into this inquiry, let's revisit some key points, including definitions and crucial results, concerning this particular class of algebras.

\begin{definition}
A \textbf{composition algebra} is an $\mathbb{F}$-algebra $\mathcal{C}$ equipped with a nondegenerate quadratic form $N: \mathcal{C} \rightarrow \mathbb{F}$ over $\mathbb{F}$ that satisfies $N(xy)=N(x)N(y)$ for any $x,y\in\mathcal{C}$. The function $N$ is referred to as the \textbf{norm} of $\mathcal{C}$.
\end{definition}

An unital composition algebra is also known as a \textbf{Hurwitz algebra}. It's important to differentiate between Hurwitz algebras that have every element invertible and those that do not. If a composition algebra contains non-zero isotropic elements (non-invertible), or equivalently, if it possesses a non-zero element with zero norm, it is termed \textbf{split}. Conversely, a composition algebra without such elements is termed \textbf{non-split}. Furthermore, it's worth noting that a Hurwitz algebra is split if and only if it is not a division algebra. An important subsequent result concerning composition algebras can be found on page 19 in \cite{springer}.

\begin{theorem}\label{key1}
There exists a unique (up to isomorphism) split Hurwitz algebra over $\mathbb{F}$ for each dimension $2$,$4$, and $8$.
\end{theorem}

Theorem \ref{key1} states that split octonion algebras are unique, therefore we can fix an algebra basis that is not dependent on the field $\mathbb{F}$.

In 2024, Bray et al. \cite{bray} constructed the following basis for the split octonion algebra. Let $I=\{\pm 0, \pm 1,\pm \omega, \pm \bar{\omega}\}$ and consider the set of symbols $B=\{e_i\}_{i \in I}$. We define multiplication on the elements of $B$ according to the following table in Figure \ref{fig:o1} and consider the eight-dimensional non-associative $\mathbb{F}$-algebra generated by $B$. In other words, we get:
\begin{itemize}
\item[i)] $e_1e_\omega = -e_\omega e_1 = e_{-\bar{\omega}}$;
\item[ii)] $e_1e_0 = e_{-0}e_1 = e_1$;
\item[iii)] $e_{-1}e_1 = -e_0$ and $e_0e_0 = e_0$;
\end{itemize}
and their images under negating all subscripts (including $0$), and multiplying all subscripts by $\omega$, where $\omega^2 = \bar{\omega}$ and $\omega\bar{\omega} = 1$. All other products of basis vectors are $0$.

\begin{figure}[h!]\label{fig2}\centering
\begin{tabular}{c|rrrrrrrr}
$  $      & $e_{-1}$ & $e_{\bar{\omega}}$    & $e_{\omega}$    & $e_0$       & $e_{-0}$    & $e_{-\omega}$       & $e_{-\bar{\omega}}$       & $e_1$ \\ \hline
$e_{-1}$ &                  $0$ & $0$    & $0$    & $0$       & $e_{-1}$    & $e_{\bar{\omega}}$       & $-e_{\omega}$       & $-e_0$ \\
$e_{\bar{\omega}}$ &                  $0$ & $0$    & $-e_{-1}$    & $e_{\bar{\omega}}$       & $0$    & $0$       & $-e_{-0}$       & $e_{-\omega}$ \\
$e_{\omega}$ &                  $0$ & $e_{-1}$    & $0$    & $e_{\omega}$       & $0$    & $-e_{-0}$       & $0$       & $-e_{-\bar{\omega}}$ \\
$e_0$ &                  $e_{-1}$ & $0$    & $0$    & $e_0$       & $0$    & $e_{-\omega}$       & $e_{-\bar{\omega}}$       & $0$ \\
$e_{-0}$ &                  $0$ & $e_{\bar{\omega}}$    & $e_{\omega}$    & $0$       & $e_{-0}$    & $0$       & $0$       & $e_1$ \\
$e_{-\omega}$ &                  $-e_{\bar{\omega}}$ & $0$    & $-e_0$    & $0$       & $e_{-\omega}$    & $0$       & $e_1$       & $0$ \\
$e_{-\bar{\omega}}$ &                  $e_{\omega}$ & $-e_0$    & $0$    & $0$       & $e_{-\bar{\omega}}$    & $-e_1$       & $0$       & $0$ \\
$e_1$ &                  $-e_{-0}$ & $-e_{-\omega}$    & $e_{-\bar{\omega}}$    & $e_1$       & $0$    & $0$       & $0$       & $0$ \\
\end{tabular}\caption{Multiplicative table of $\mathbb{O}_{\mathfrak{s}}$}\label{fig:o1}\end{figure}

Note that the multiplicative identity is given by $e_0 + e_{-0} = 1$. For an element $x= \sum_{i \in I} \lambda_i e_i$ of this algebra we define a trace:
$$
T(x)=\lambda_0+ \lambda_{-0}
$$
and norm:
$$
N(x)=\lambda_{-1}\lambda_{1}+ \lambda_{-\bar{\omega}}\lambda_{\bar{\omega}}+ \lambda_{-\omega}\lambda_{\omega}+ \lambda_{-0}\lambda_{0}.
$$
Then it is easy to see that:
\begin{equation*}
x^2=T(x)x-N(x).
\end{equation*}
Thus, for $x$ such that $N(x)\neq0$, we have:
\begin{equation*}
x^{-1}=N(x)^{-1}(T(x)-x).
\end{equation*}

Now we will state the main answer to the Conjecture \ref{c1}.

\begin{theorem}\label{t2}
Let $D$ be a split octonion algebra over a field $\mathbb{F}$, and suppose $f, g : D \rightarrow D$ are additive maps satisfying the identity \eqref{vukmana}:
\begin{equation*}
f(x) + x^2g(x^{-1}) = 0
\end{equation*}
for all $x \in D^\times$. Then there exists a fixed element $q \in D$ such that $f(x) = xq$ and $g(x) = -xq$ for all $x \in D$.
\end{theorem}

We will divide it into some cases. First we consider the case where $\mathrm{Char}(\mathbb{F})\neq2$.

\begin{proof}
Let $h(x)=f(x)-xf(1)$. Then $h$ is additive, $h(1)=0$ and, by Step 1 of the proof of Theorem \ref{t1}, for all $x\in D$ such that $x\in D^\times$ and $1-x\in D^\times$ we have:
\begin{equation}
\label{4a}
h(x^2)=2xh(x).
\end{equation}
Now we divide it into some steps.

\textbf{Step 1:} We will prove that $h(\alpha e_i)=0$ for every $\alpha\in\mathbb{F}$ and $i\in\{\pm1,\pm\omega,\pm\bar{\omega}\}$. Indeed, let $\alpha\in\mathbb{F}$, $i\in\{\pm1,\pm\omega,\pm\bar{\omega}\}$ and:
\begin{equation*}
x=e_0-e_{-0}+e_{i\omega}-e_{-i\omega}+\alpha e_i.
\end{equation*}
Then:
\begin{equation*}
1-x=2e_{-0}-e_{i\omega}+e_{-i\omega}-\alpha e_i.
\end{equation*}
Thus $T(x)=0$, $N(x)=-2$ and $N(1-x)=-1$, so, by $\mathrm{Char}(\mathbb{F})\neq2$, \eqref{4a} holds. But:
\begin{equation*}
x^2=T(x)x-N(x)=2.
\end{equation*}
Hence by \eqref{4a} we have $h(2)=2xh(x)$. But $h(1)=0$, $\mathrm{Char}(\mathbb{F})\neq2$ and $x\in D^\times$, so $h(x)=0$, in other words:
\begin{equation}
\label{4b}
h(e_0-e_{-0}+e_{i\omega}-e_{-i\omega}+\alpha e_i)=0.
\end{equation}
The above equation holds for all $\alpha\in\mathbb{F}$, in particular:
\begin{equation}
\label{4b1}
h(e_0-e_{-0}+e_{i\omega}-e_{-i\omega})=0.
\end{equation}
Therefore, for any $\alpha\in\mathbb{F}$, we have \eqref{4b} and \eqref{4b1}, so we conclude that $h(\alpha e_i)=0$.

\textbf{Step 2:} We will prove that $h(e_{i})=0$ for $i\in\{\pm0\}$. Indeed, by \eqref{4b1} and Step 1 we obtain $h(e_0)-h(e_{-0})=0$. Also $h(1)=0$, so we have $h(e_0)+h(e_{-0})=0$. Because $\mathrm{Char}(\mathbb{F})\neq2$, we obtain $h(e_0)=h(e_{-0})=0$.

\textbf{Step 3:} We will prove that $h(\alpha)=0$ for any $\alpha\in\mathbb{F}$. Indeed $h(0)=0$, and let $\alpha\in\mathbb{F}\setminus\{0\}$ and:
\begin{equation*}
x=e_0+\alpha e_1+e_{-1}.
\end{equation*}
Then:
\begin{equation*}
1-x=e_{-0}-\alpha e_1-e_{-1}.
\end{equation*}
Thus $T(x)=1$ and $N(x)=N(1-x)=\alpha$, so \eqref{4a} holds. But:
\begin{equation*}
x^2=T(x)x-N(x)=x-\alpha.
\end{equation*}
Hence, by \eqref{4a}, we have $h(x-\alpha)=2xh(x)$, but by Steps 1 and 2 we have $h(x)=0$, so we obtain $h(\alpha)=0$.

\textbf{Step 4:} We will prove that $h(\alpha e_i)=0$ for all $\alpha\in\mathbb{F}$ and $i\in\{\pm0\}$. Indeed let $\alpha\in\mathbb{F}$ and:
\begin{equation*}
x=\alpha e_0-\alpha e_{-0}+(\alpha^2+1)e_1+e_{-1}.
\end{equation*}
Then:
\begin{equation*}
1-x=-(\alpha-1)e_0+(\alpha+1)e_{-0}-(\alpha^2+1)e_1-e_1.
\end{equation*}
Thus $T(x)=0$, $N(x)=1$ and $N(1-x)=2$, so, by $\mathrm{Char}(\mathbb{F})\neq2$, \eqref{4a} holds. But:
\begin{equation*}
x^2=T(x)x-N(x)=-1.
\end{equation*}
Hence, by \eqref{4a}, we have $h(-1)=2xh(x)$. But $h(1)=0$, $\mathrm{Char}(\mathbb{F})\neq2$ and $x\in D^\times$, so $h(x)=0$. Thus by Step 1 we obtain $h(\alpha e_0)-h(\alpha e_{-0})=0$. Also by Step 3 we have $h(\alpha e_0)+h(\alpha e_{-0})=0$. Because $\mathrm{Char}(\mathbb{F})\neq2$, we obtain $h(\alpha e_0)=h(\alpha e_{-0})=0$.

\textbf{Step 5:} We will prove that $f(x)=xf(1)$ for all $x\in D$. Indeed, for every element $x\in D$ is of the form $x=\sum_{i\in I}\lambda_ie_i$, so:
\begin{equation*}
h(x)=h\left(\sum_{i\in I}\lambda_ie_i\right)=\sum_{i\in I}h(\lambda_ie_i)=0.
\end{equation*}
In other words, for every $x\in D$ we have $h(x)=0$, so $f(x)=xf(1)$.

\textbf{Step 6:} We will conclude the proof. Indeed, for every $x\in D^\times$, by \eqref{vukmana} with $x^{-1}$ rather than $x$, we obtain:
\begin{equation*}
f(x^{-1})+x^{-2}g(x)=0,
\end{equation*}
we get:
\begin{equation*}
g(x)+x^2f(x^{-1})=0.
\end{equation*}
This is \eqref{vukmana} with $f$ and $g$ interchanged, so we can prove a result for $g$ analogous to Step 5, namely that $g(x)=xg(1)$ for all $x\in D$. By \eqref{vukmana} with $x=1$, we have $f(1)+g(1)=0$, so that $g(1)=-f(1)$, therefore we can let $q=f(1)$.
\end{proof}

Now we consider the case $\mathrm{Char}(\mathbb{F})=2$.

\begin{proof}
Let $h(x)=f(x)+xf(1)$ and $k(x)=g(x)+xg(1)$. Then the functions $h,k$ are additive and, by \eqref{vukmana} with $x=1$, so that $f(1)+g(1)=0$, the functions $h,k$ satisfy the equation:
\begin{equation}
\label{vukman1}
h(x)+x^2k(x^{-1})=0
\end{equation}
for all $x\in D^\times$. Also $h(1)=k(1)=0$ and, by Step 1 of the proof of Theorem \ref{t1}, for all $x\in D$ such that $x\in D^\times$ and $x+1\in D^\times$ we have:
\begin{equation}
\label{4a'}
h(x^2)=0.
\end{equation}
Now we divide it into some steps.

\textbf{Step 1:} We will prove that $h(\alpha e_i)=0$ for every $\alpha\in\mathbb{F}$ and $i\in\{\pm1,\pm\omega,\pm\bar{\omega}\}$. Indeed, let $\alpha\in\mathbb{F}$, $i\in\{\pm1,\pm\omega,\pm\bar{\omega}\}$ and:
\begin{equation*}
x=e_0+e_{i\omega}+e_{-i\omega}+\alpha e_i.
\end{equation*}
Then:
\begin{equation*}
x+1=e_{-0}+e_{i\omega}+e_{-i\omega}+\alpha e_i.
\end{equation*}
Thus $T(x)=1$ and $N(x)=N(x+1)=1$, so \eqref{4a'} holds. But:
\begin{equation*}
x^2=T(x)x+N(x)=x+1.
\end{equation*}
Hence, by \eqref{4a'}, we have $h(x+1)=0$, but $h(1)=0$, so $h(x)=0$, in other words:
\begin{equation}
\label{4b'}
h(e_0+e_{i\omega}+e_{-i\omega}+\alpha e_i)=0.
\end{equation}
The above equation holds for any $\alpha\in\mathbb{F}$, in particular:
\begin{equation}
\label{4b1'}
h(e_0+e_{i\omega}+e_{-i\omega})=0.
\end{equation}
Therefore, for any $\alpha\in\mathbb{F}$, we have \eqref{4b'} and \eqref{4b1'}, so $h(\alpha e_i)=0$.

\textbf{Step 2:} We will prove that $h(e_i)=0$ for $i\in\{\pm0\}$. Indeed, by \eqref{4b1'} and Step 1 we have $h(e_0)=0$. Also $h(1)=0$, so we have $h(e_0)+h(e_{-0})=0$. Therefore $h(e_{-0})=0$.

\textbf{Step 3:} We will prove that $h(\alpha)=0$ for all $\alpha\in\mathbb{F}$. Indeed $h(0)=0$, and let $\alpha\in\mathbb{F}\setminus\{0\}$ and:
\begin{equation*}
x=e_0+\alpha e_1+e_{-1}.
\end{equation*}
Then:
\begin{equation*}
x+1=e_{-0}+\alpha e_1+e_{-1}.
\end{equation*}
Thus $T(x)=1$ and $N(x)=N(x+1)=\alpha$, so \eqref{4a'} holds. But:
\begin{equation*}
x^2=T(x)x+N(x)=x+\alpha.
\end{equation*}
Hence, by \eqref{4a'} we have $h(x+\alpha)=0$, but by Steps 1 and 2, we have $h(x)=0$, so we obtain $h(\alpha)=0$.

\textbf{Step 4:} For every $x\in D^\times$, then by \eqref{vukman1} with $x^{-1}$ in place of $x$ we have:
\begin{equation*}
h(x^{-1})+x^{-2}k(x)=0,
\end{equation*}
so we obtain:
\begin{equation*}
k(x)+x^2h(x^{-1})=0.
\end{equation*}
This is \eqref{vukman1} with $h$ and $k$ interchanged, so we can prove results for $k$ that are analogous to Steps 1 to 3.

\textbf{Step 5:} We will prove that $h(\alpha e_i)=k(\alpha e_i)=0$ for all $\alpha\in\mathbb{F}$ and $i\in\{\pm0\}$. By Steps 2 and 4, we have $h(\alpha e_i)=k(\alpha e_i)=0$ for $\alpha\in\{0,1\}$ and $i\in\{\pm0\}$. Now let $\alpha\in\mathbb{F}\setminus\{0,1\}$, $\beta\in\mathbb{F}\setminus\{0\}$, $i\in\{\pm0\}$ and:
\begin{equation*}
x=\alpha e_i+\beta e_1+\beta^{-1} e_{-1}.
\end{equation*}
Then $T(x)=\alpha$ and $N(x)=1$, so $x\in D^\times$ and:
\begin{equation*}
x^{-1}=N(x)^{-1}(T(x)+x)=\alpha e_{-i}+\beta e_1+\beta^{-1} e_{-1},
\end{equation*}
\begin{equation*}
x^2=T(x)x+N(x)=(\alpha^2+1)e_i+e_{-i}+\alpha\beta e_1+\alpha\beta^{-1}e_{-1}.
\end{equation*}
By \eqref{vukman1}, we have $h(x)=x^2k(x^{-1})$, but by Steps 1 and 4 we have $h(x)=h(\alpha e_i)$ and $k(x^{-1})=k(\alpha e_{-i})$, so:
\begin{equation}
\label{z}
h(\alpha e_i)=((\alpha^2+1)e_i+e_{-i}+\alpha\beta e_1+\alpha\beta^{-1}e_{-1})k(\alpha e_{-i}).
\end{equation}
In particular, applying the \eqref{z} for $\beta=\alpha$ and for $\beta=\alpha^{-1}$, we obtain:
\begin{equation*}
h(\alpha e_i)=((\alpha^2+1)e_i+e_{-i}+\alpha^2 e_1+e_{-1})k(\alpha e_{-i}),
\end{equation*}
\begin{equation*}
h(\alpha e_i)=((\alpha^2+1)e_i+e_{-i}+ e_1+\alpha^2e_{-1})k(\alpha e_{-i}).
\end{equation*}
By summing the above two equations, we obtain:
\begin{equation*}
0=((\alpha^2+1)e_1+(\alpha^2+1)e_{-1})k(\alpha e_{-i})
\end{equation*}
But:
\begin{equation*}
N((\alpha^2+1)e_1+(\alpha^2+1)e_{-1})=(\alpha^2+1)^2=(\alpha+1)^4\neq0.
\end{equation*}
Therefore $k(\alpha e_{-i})=0$. Now by \eqref{z} we have $h(\alpha e_i)=0$.

\textbf{Step 6:} We will conclude the proof. Let $l\in\{h,k\}$. For every element $x\in D$, then $x$ is of the form $x=\sum_{i\in I}\lambda_ie_i$, so:
\begin{equation*}
l(x)=l\left(\sum_{i\in I}\lambda_ie_i\right)=\sum_{i\in I}l(\lambda_ie_i)=0.
\end{equation*}
In other words, for every $x\in D$ we have $l(x)=0$, so $f(x)=xf(1)$ and $g(x)=xg(1)$
for all $x\in D$. Now take $q=f(1)$ and remember that $f(1)+g(1)=0$.
\end{proof}

\end{document}